\documentclass[11pt]{amsart}

\usepackage{amsmath,amsthm, amssymb}
\usepackage[T1]{fontenc}
\usepackage{latexsym}
\usepackage{subfigure, epsfig, epsf}
\usepackage{color,fancyhdr,lastpage,graphicx}
\usepackage{xspace}
\usepackage{wasysym}
\usepackage{setspace}
\usepackage{bbm}
\usepackage{stmaryrd}

\usepackage{pstricks}
\usepackage{pst-node}
\usepackage{pst-tree}

\newcommand{\RR}{\mathbb{R}}    \newcommand{\LL}{\mathbb{L}}

      \newcommand{\XX}{\mathbb{X}}   \newcommand{\bbA}{\mathbb{A}}
\newcommand{\bbW}{\mathbb{W}}

\newcommand{\mcP}{\mathcal{P}} \newcommand{\mcC}{\mathcal{C}}  \newcommand{\mcF}{\mathcal{F}}

  \newcommand{\be}{{\bf e}}

    \newcommand{\bfX}{{\bf X}}

\newcommand{\bbX}{\mathbb{X}}

 \newcommand{\PP}{\mathbb{P}}

\newcommand{\ep}{\epsilon}
\renewcommand{\leq}{\leqslant}
\renewcommand{\geq}{\geqslant}

\newcommand{\wrt}{with respect to }

\newcommand{\ssk}{\smallskip}

\newtheorem{thm}{\hspace{-0.15cm}  {\sc Theorem} }
\newtheorem{cor}[thm]{\hspace{-0.15cm}  {\sc Corollary} }
\newtheorem{prop}[thm]{\hspace{-0.15cm} {\sc Proposition}}
\newtheorem{defn}[thm]{ \hspace{-0.3cm} {\sc Definition}}

\newtheorem{rem}[thm]{{\sc Remark}}
\newtheorem{rems}[thm]{{\sc Remarks}}

\numberwithin{equation}{section} 

\newenvironment{Dem}{%
    \begin{list}{\hspace{0.5cm}{\sc Proof --}}{%
        \setlength{\topsep}{0pt}%
        \setlength{\leftmargin}{0pt}%
        \setlength{\rightmargin}{0pt}%
        \setlength{\listparindent}{0pt}%
        \setlength{\itemindent}{0pt}%
        \setlength{\parsep}{0pt}%
        \addtolength{\leftmargin}{20pt}%
        \addtolength{\rightmargin}{0pt}%
    } \item }{\hfill{\space $\rhd$}\end{list}\smallskip}

\title{Path-dependent rough differential equations}
\date{\today}
\author{I. Bailleul}
\address{IRMAR, 263 Avenue du General Leclerc, 35042 RENNES, France}
\email{ismael.bailleul@univ-rennes1.fr}

\begin{document}

\maketitle

\begin{abstract}
We show in this work how the machinery of $\mcC^1$-approximate flows introduced in \cite{RPDrivenFlows} provides a very efficient tool for proving well-posedness results for path-dependent rough differential equations on flows of the form 
$$
d\varphi = Vh(dt) + \textrm{F}{\bfX}(dt),
$$
for smooth enough path-dependent vector fields $V,\textrm{F}=\big(V_1,\dots,V_\ell\big)$, any H\"older weak geometric $p$-rough path $\bfX$ and any $\alpha$-H\"older path $h$, with $\alpha+\frac{1}{p}>1$.
\end{abstract}

\section{Introduction}
\label{SectionIntroduction}

Although it makes no doubt that the pointwise Markovian description of numerous physical systems given by ordinary differential equations on some state space provides an excellent approximate description of its evolution, many physical systems have characteristics whose internal states evolve along their history and influence their response to their environment. The monographs \cite{Mohammed, Mohammed98} provides examples ranging from engineering and physics to biology. The seminal works of Ito-Nisio \cite{ItoNisio} and Kushner \cite{Kushner} on path-dependent stochastic differential equations took a mature form in the work of Mohammed \cite{Mohammed} and others on stochastic functional differential equations. However, the classical setting of Ito integration used in these works appears as too narrow to handle a number of problems coming from numerous models in applied mathematics, involving for instance integrals with respect to Gaussian processes such as fractional Brownian motion.

The conceptual framework and tools developped by T. Lyons in his theory of rough paths \cite{Lyons97} offer a flexible and convenient setting to give some meaning and study differential equations of the form
\begin{equation}
\label{EqGenericPDDE}
dx_t = V\big(x_{[0,t]}\big)dt + \textrm{F}\big(x_{[0,t]}\big){\bfX}(dt),
\end{equation}
for some sufficiently regular path-dependent maps $V$,F and a driving signal $\bfX$ which can be far rougher than Brownian motion. Several approaches to rough paths and rough differential equations are now available. We shall use in this work the machinery of $\mcC^1$-approximate flows, developped in \cite{RPDrivenFlows} to recover and extend most of the basic results on rough differential equations in a Banach setting. This will enable us to deal with systems with bounded and unbounded memory in an equal footing, proving well-posedness results for equation \eqref{EqGenericPDDE} under sufficient regularity conditions on the path-dependent maps $V$ and $F$. 

\ssk

After recalling in section \ref{SectionApproximateFlows} the results on flows and $\mcC^1$-approximate flows proved in \cite{RPDrivenFlows}, we dedicate section \ref{SectionPathDependentRDE} to investigating the well-posedness problem for the equation
\begin{equation}
\label{EqGeneralPDRDE}
d\varphi = Vh(dt) + \textrm{F}{\bfX}(dt)
\end{equation}
on flows on $\RR^d$, where $h : [0,T]\rightarrow\RR$ is an $\alpha$-Lipschitz real-valued path and ${\bfX}=\big(X,\bbX\big)$ is a weak geometric $p$-rough path, with $2\leq p<3$ and $\alpha+\frac{1}{p}>1$. The set of path-dependent vector fields considered in this work, and their regularity properties, are defined in section \ref{SubsectionPathDpdtVectorFields}. We show in section \ref{SubsectionPathDpdtRDEs} how one can define a natural $\mcC^1$-approximate flow from $V$, F and $\bfX$, which leads, together with the results of section \ref{SectionApproximateFlows}, to the well-posed character of equation \eqref{EqGeneralPDRDE}. This result comes with a number of corollaries: a convergence theorem for sequences equations with vanishing delay, sharp pathwise convergence rates for Milstein-types schemes, and large deviation results and support theorem for stochastic path-dependent differential equations, extending earlier results.

We extend these results to any weak geometric $p$-rough path, under the condition $\alpha+\frac{1}{p}>1$, in section \ref{SectionGeneralRoughPaths}. 

\medskip

\noindent {\bf Notations.} We gather in the remainder of this introduction some notations which are used throughout the article.

\medskip

$\bullet$ We denote by $\|\cdot\|_\infty$the supremum norm over some fixed time interval $[0,T]$. We also denote by $c$ a constant whose value may change from place to place.

\ssk

$\bullet$ Denote by $(e_1,\dots,e_\ell)$ the canonical basis of $\RR^\ell$, seen as a subset of the tensor algebra $T^{(\infty)}_\ell = \bigoplus_{n\geq 1} \big(\RR^\ell\big)^{\otimes n}$ over $\RR^\ell$ -- see for instance \cite{Lyons97}. For a tuple $I=(i_1,\dots,i_r)$, write
$$
{\be}_{[I]}=\Big[e_{i_1},\big[e_{i_2},\dots[e_{i_{r-1}},e_{i_r}]\big]\dots\Big],\quad \textrm{ and   } {\be}_I=e_{i_1}e_{i_2}\cdots e_{i_r},
$$
where the above products are in $T^{(\infty)}_\ell$. Write $T_\ell^{[p]}$ for the truncated tensor algebra of depth $[p]$, and $Y^I$ for the coordinates of an element $Y$ of $T_\ell^{[p]}$ in its canonical basis.

\ssk

$\bullet$ Given a Holder weak geometric $p$-rough path $\bfX$, living in the step-$[p]$ free nilpotent Lie group over $\RR^\ell$, denote by ${\bf \Lambda}_{ts} = 0\oplus \Lambda^1_{ts}\oplus\cdots\oplus \Lambda^{[p]}_{ts}$, the logarithm in $T^{[p]}_\ell$ of ${\bfX}_{ts}$
\begin{equation}
\label{EqDefnLambda}
\exp\Big({\bf \Lambda}_{ts}\Big) = X^I_{ts} {\be}_I = \bfX_{ts},
\end{equation}
for all $0\leq s\leq t\leq T$; it takes values in the finite dimensional step-$[p]$ free nilpotent Lie algebra over $\RR^\ell$. Notice that since ${\bf \Lambda}$ is polynomial in $\bfX$, it is a continuous function $\bfX$.


\section{Flows and $\mcC^1$-approximate flows}
\label{SectionApproximateFlows}

Following earlier works of Lyons \cite{Lyons97} and Gubinelli \cite{Gubinelli}, Feyel and de la Pradelle introduced in \cite{FdlP1} a notion of almost-additive function, with values in a Banach space $E$. These are $E$-valued functions $\mu$ defined on the simplex $\big\{(s,t)\in [0,T]^2\,;\,s\leq t\big\}$, and such that we have
$$
\Big|\big(\mu_{tu}+\mu_{us}\big)-\mu_{ts}\Big| \leq c|t-s|^\gamma,
$$
for some positive constants $c$ and $\gamma>1$, and all $0\leq s\leq u\leq t\leq T$. A function $\mu$ for which the above left hand side is identically null gives the increments of another $E$-valued function $\varphi$ defined on the interval $[0,T]$, in the sense that $\mu_{ts} = \varphi_t-\varphi_s$, for all $0\leq s\leq t\leq T$. Feyel and de la Pradelle show that one can associate to any almost-additive function $\mu$ a unique function $\varphi : [0,T]\rightarrow E$, whose increments are close to $\mu$, in the sense that the inequality
$$
\big|\big(\varphi_t-\varphi_s\big)-\mu_{ts}\big| \leq c|t-s|^\gamma
$$
holds for some positive constant $c$, for all $0\leq s\leq t\leq T$. Thinking for instance of $\mu_{ts}$ as $f_s\big(g_t-g_s\big)$, for two real-valued functions $f$ and $g$, respectively $\alpha$ and $\beta$-H\"older continuous, with $\alpha+\beta>1$, the function $\varphi$ given by their theorem is the well-defined Young integral $\int_0^\bullet f_sdg_s$; see \cite{LyonsStFlour}.

Feyel, de la Pradelle and Mokobodzki gave in \cite{FdlP2} a multiplicative extension of that result, well-adapted to Banach algebra-valued maps $\mu$, but not flexible enough to construct some large classes of flows. The notion of $\mcC^1$-approximate flow introduced in \cite{RPDrivenFlows} provides a convenient and flexible different extension of their result, taylor made to construct flows of maps on Banach spaces. We give here a simplified definition of a $\mcC^1$-approximate flow sufficient for our needs in that work. Recall that a {\it flow} is a family $\big(\varphi_{ts}\big)_{0\leq s\leq t\leq T}$ of maps from $\RR^d$ to itself such that one has $\varphi_{tu}\circ\varphi_{us} = \varphi_{ts}$, for all $s\leq u\leq t$.

\begin{defn}
\label{DefnC1ApproximateFlow}
A \emph{{\bf $\mcC^1$-approximate flow}} is a family $\big(\mu_{ts}\big)_{0\leq s\leq t\leq T}$ of $\mcC^2$ maps from $\RR^d$ to itself, depending continuously on $(s,t)$ in the topology of uniform convergence, such that 
   \begin{equation}
   \label{EqRegularityBounds}
   \big\|\mu_{ts} - \textrm{\emph{Id}}\big\|_{\mcC^2} \leq o_{t-s}(1),
   \end{equation} 
and there exists some positive constants $c_1$ and $\gamma>1$, such that 
\begin{equation}
\label{EqMuMu}
\big\|\mu_{tu}\circ\mu_{us}-\mu_{ts}\big\|_{\mcC^1} \leq c_1\,|t-s|^\gamma
\end{equation}
hold for any $0\leq s\leq u\leq t\leq T$.
\end{defn}

Given a partition $\pi_{ts}=\{s=s_0<s_1<\cdots<s_{n-1}<s_n=t\}$ of $(s,t)\subset [0,T]$, set 
$$
\mu_{\pi_{ts}} = \mu_{s_ns_{n-1}}\circ\cdots\circ\mu_{s_1s_0}.
$$
The following statement, proved in \cite{RPDrivenFlows}, provides a flexible tool for constructing flows from $\mcC^1$-approximate flows, while generalizing the above mentionned result of Feyel-de la Pradelle to the non-commutative setting of $\mcC^2$ maps on $\RR^d$.
 It holds more generally on any Banach space under weaker conditions on $\mu$ than in definition \ref{DefnC1ApproximateFlow}.
 
\begin{thm}[Constructing flows on $\RR^d$]
\label{ThmConstructingFlows}
A $\mcC^1$-approximate flow defines a unique flow $(\varphi_{ts})_{0\leq s\leq t\leq T}$ on $\RR^d$ such that 
$$
\big\|\varphi_{ts}-\mu_{ts}\big\|_\infty \leq c|t-s|^\gamma
$$
holds for some positive constant $c$, for all $0\leq s\leq t\leq T$ sufficiently close, say $|t-s|\leq \delta$; this flow satisfies the inequality 
\begin{equation}
\label{EqApproxVarphiMu}
\big\|\varphi_{ts}-\mu_{\pi_{ts}}\big\|_\infty\leq c_1^2 \,T\,\big|\pi_{ts}\big|^{\gamma-1}
\end{equation} 
for any partition $\pi_{ts}$ of any interval $(s,t)\subset [0,T]$, with mesh $\big|\pi_{ts}\big|\leq \delta$. 
\end{thm}

This result was used in \cite{RPDrivenFlows} to recover and extend in a simple way most of the basic results on rough differential equations driven by weak geometric rough paths: a well-posedness result for rough differential equations on flows was proved in a Banach setting, giving back and extending some results for classical rough differential equations, and leading to a Peano existence theorem under minimal regularity conditions on the driving vector fields. Higher order Euler expansions and non-explosion under linear growth conditions on the vector fields were also proved; an application to some stochastic rough differential equation of mean field type was also given.

\section{Path-dependent rough differential equations}
\label{SectionPathDependentRDE}

We show in this section how the machinery of $\mcC^1$-approximate flows recalled in section \ref{SectionApproximateFlows} can be used to give some meaning and solve some path-dependent rough differential equations of the form
\begin{equation}
\label{EqPathDpdtRDE}
d\varphi = Vh(dt) + \textrm{F}{\bfX}(dt),
\end{equation}
where $h$ is an $\alpha$-H\"older real-valued path, $\bfX$ is a H\"older weak geometric $p$-rough path, with $2<p<3$ and $\alpha+\frac{1}{p}>1$, and $V$ and $\textrm{F} = \big(V_1,\dots,V_\ell\big)$ are path-dependent vector fields over $\RR^d$. Our approach consists, as in \cite{RPDrivenFlows}, in associating to $(s,t),h,\bfX$ and $V$, F, an ordinary (path-dependent) differential equation whose time $1$ map $\mu_{ts}$ will be shown to define a $\mcC^1$-approximate flow with the awaited Euler expansion. The unique flow associated to $\mu$ by theorem \ref{ThmConstructingFlows} will be said to be the solution to the path-dependent rough differential equation \eqref{EqPathDpdtRDE}. Some regularity conditions on the path-dependent vector fields $V,V_1,\dots,V_\ell$ are needed to make that approach work. We define them in section \ref{SubsectionPathDpdtVectorFields}, introducing a new notion of fine-differentiability. The above program for solving equation \eqref{EqPathDpdtRDE} is run in section \ref{SubsectionPathDpdtRDEs}.

\subsection{Path-dependent vector fields on $\RR^d$}
\label{SubsectionPathDpdtVectorFields}

Set 
$$
\mcP := \bigcup_{0\leq a\leq T}\mcC\Big([0,a],\RR^d\Big),
$$
with $\mcC\big([0,0],\RR^d\big)$ identified to $\RR^d$; denote by $y_{[0,a]}$ a generic element of $\mcC\big([0,a],\RR^d\big)$. Given $y_{[0,a]}\in\mcP$, write $\overline{y}_\bullet$ for the path over $[0,T]$ whose values coincide with $y_{[0,a]}$ over $[0,a]$, and which is constant, equal to  $y_a$, on the interval $[a,T]$. Note that for any given $y_{[0,a]}$ and $z_{[0,a]}$, we have $\sup_{0\leq t\leq a}\;|y_t-z_t| = \big\|\overline{y}_\bullet-\overline{z}_\bullet\big\|_\infty$. We define a metric d on $\mcP$ setting 
$$
\textrm{d}\Big(y_{[0,a]},z_{[0,b]}\Big) = |b-a| + \big\|\overline{y}_\bullet - \overline{z}_\bullet\big\|_\infty.
$$
Given $x\in\RR^d$, we denote by $\overline{x}_{[0,a]}$ the constant path over the time interval $[0,a]$, identically equal to $x$.

\begin{defn}
\label{DefnPathDpdtVectorField}
A {\bf path-dependent vector field on $\RR^d$} is a map $V$ from $\mcP$ to $\RR^d$. 
\begin{itemize}
   \item It is said to be {\bf Lipschitz} if it satisfies the inequality
$$
\Big|V\Big(y_{[0,a]}\Big)-V\Big(z_{[0,a]}\Big)\Big| \leq \lambda\|\overline{y}_\bullet-\overline{z}_\bullet\|_\infty,
$$
for some positive constant $\lambda$, for all $y_{[0,a]}$ and $z_{[0,a]}$ in $\mcP$, and any $0\leq a\leq T$. It is said to be {\bf strongly Lipschitz} if we have
$$
\Big|V\Big(y_{[0,a]}\Big)-V\Big(z_{[0,b]}\Big)\Big| \leq \lambda\|\overline{y}_\bullet-\overline{z}_\bullet\|_\infty,
$$
for some positive constant and all $y_{[0,a]}$ and $z_{[0,b]}$ in $\mcP$, with $0\leq a,b\leq T$.
   \item It is said to be {\bf of class $\mcC^k$} if its restriction to each space $\mcC\big([0,a],\RR^d\big)$ is $\mcC^k$-Fr\'echet differentiable, with respect to the supremum norm on $[0,a]$.
\end{itemize}
\end{defn}

\noindent If $V$ is bounded and $\mcC^k$, for some $k\geq 1$, with bounded Fr\'echet derivatives $F^{(1)},\dots,F^{(k)}$, we define its norm to be 
$$
\|F\|_{\mcC^k} = \sum_{i=0}^k \Big\|F^{(i)}\Big\|,
$$
with $F^{(0)}$ standing for $F$ itself. We denote by $\mcC^k_b$ the set of maps for which $\|F\|_{\mcC^k}$ is  finite. As expected, a path-dependent vector field $V$ is Lipschitz if $V$ is $\mcC^1$, with a bounded derivative. The restriction to $\mcC\big([0,0],\RR^d\big)$ of a $\mcC^k$ path-dependent vector field is naturally identified with a $\mcC^k$ vector field on $\RR^d$. The following examples of path-dependent vector fields are commonly encountered in applications:
\begin{equation}
\label{EqExample1}
V\Big(y_{[0,a]}\Big) = f\Big(y_{(a-r_1)\vee 0},\dots,y_{(a-r_n)\vee 0}\Big),
\end{equation}
for some $\RR^d$-valued function $f$ on $\big(\RR^d\big)^n$, for $n\geq 1$, and $r_1,\dots,r_n$ either non-negative constants or non-negative functions of $a$;
\begin{equation}
\label{EqExample2}
V\Big(y_{[0,a]}\Big) = \frac{1}{a}\int_0^a g\big(y_r\big)dr,
\end{equation}
for $a>0$, and $V\Big(y_{[0,0]}\Big) = g(y_0)$, for some $\RR^d$-valued function $g$ on $\RR^d$. Both vector fields are of class $\mcC^k$ if $f$ and $g$ are $\mcC^k$. On the other hand, the path-dependent vector field 
$$
V\Big(y_{[0,a]}\Big) = \int_0^{(a-r)\vee 0} x_{a-s}\,dW_s
$$
corresponding to stochastic white noise perturbation of the memory does not fit into the framework of definition \ref{DefnPathDpdtVectorField}.

\ssk

Given a path $y_{[0,a]}\in\mcP$ and a path-dependent vector field $V$, set 
$$
y^{V,\ep}_{[0,a]} : s\in [0,a]\mapsto y_s + \ep V\Big(y_{[0,a]}\Big).
$$
For a $\mcC^1$ path-dependent vector field $W$, set 
\begin{equation}
\label{EqDefnVW}
(VW)\big(y_{[0,a]}\big) = \frac{d}{d\ep}_{\big|\ep=0} W\Big(y^{V,\ep}_{[0,a]}\Big);
\end{equation}
this definition coincides with the usual differentiation when $a=0$ and the path $y_{[0,0]}$ is identified with the point $y_0$ of $\RR^d$.  We adopt this definition of $VW$ rather than the traditional
\begin{equation}
\label{EqOtherDefnVW}
\frac{d}{d\ep}_{\big|\ep=0} W\left(\Big\{y_\bullet+\ep V\big(y_{[0,\bullet]}\big)\Big\}_{[0,a]}\right),
\end{equation}
since the latter depends on the whole path $V\big(y_{[0,\bullet]}\big)_{[0,a]}$ rather than just on $V\big(y_{[0,a]}\big)$, as in \eqref{EqDefnVW}; definition \eqref{EqOtherDefnVW} leads to different dynamics for which all the claims below hold true.

\ssk

Write $D_\bullet W$ for the linear map $V\mapsto VW$. We define the bracket of $V$ and $W$ to be 
$$
[V,W]\big(y_{[0,a]}\big) = (VW)\Big(y_{[0,a]}\Big) - (WV)\Big(y_{[0,a]}\Big),
$$
if $V$ is also a $\mcC^1$ path-dependent vector field; this definition is coherent with the usual definition of the bracket for $a=0$.

\ssk

\noindent It is well-known that if $V$ is Lipschitz then the path-dependent ordinary differential equation 
\begin{equation}
\label{EqPathDpdtODE}
y_t = y_0 + \int_0^tV\Big(y_{[0,s]}\Big)ds
\end{equation}
on the time interval $[0,T]$, has a unique solution which defines a $\mcC^k$ function of $y_0$ if $V$ is of class $\mcC^k$. As in the classical setting, if $V$ depends continuously on some parameter $\lambda$, in the sense that the $\mcC\big([0,T],\RR^d\big)$-valued map 
$$
\Big(z_{[0,T]},\lambda\Big) \mapsto y_0+ \int_0^\bullet V_\lambda\Big(z_{[0,s]}\Big)ds
$$
is continuous, then the solution $y_\bullet$ to equation \eqref{EqPathDpdtODE} depends continuously on $\lambda$ in $\|\cdot\|_\infty$-norm.

\medskip

\noindent Set 
$$
\mcP^1 := \bigcup_{0\leq a\leq T}\mcC^1\Big([0,a],\RR^d\Big),
$$
with $\mcC^1\big([0,0],\RR^d\big)$ identified with $\RR^d\times\RR^d$, and write $\dot y_\bullet$ for the derivative of an $\RR^d$-valued $\mcC^1$ path $y_\bullet$.

Given a Banach space $E$, denote by $\textrm{L}(\mcP,E)$ the set of continuous linear maps from $\mcP$ to $E$, and identify $\textrm{L}\big(\mcP,\textrm{L}(\mcP,E)\big)$ to $\textrm{L}\big(\mcP^{\otimes 2},E\big)$, with a similar meaning for $\textrm{L}\big(\mcP^{\otimes n},E\big)$, for $n\geq 2$. 

\begin{defn}
A map $F$ from $\mcP$ to some Banach space $E$ is said to be {\bf finely-$\mcC^1$} if there exists a continuous map $F' : \mcP\mapsto\textrm{\emph{L}}(\mcP,E)$ such that
\begin{equation}
\label{EqFineRegularity}
F\big(y_{[0,a]}\big) = F\big(y_{[0,0]}\big) + \int_0^a F'\Big(y_{[0,r]};\dot y_{[0,r]}\Big)dr,
\end{equation}
holds for any $y_{[0,a]}\in\mcP^1$. Write $F^{\{1\}}$ for its fine derivative $F'$. We say that $F$ is {\bf finely}-$\mcC^k$, with $k\geq 2$, if $F^{\{k-1\}}$ is finely-$\mcC^1$. 
\end{defn}

If $F$ is bounded and finely-$\mcC^k$, for some $k\geq 1$, with bounded derivatives $F^{\{1\}},\dots,F^{\{k\}}$, we define its norm to be 
$$
\big|F\big|_{\mcC^k} = \sum_{i=0}^k \Big|F^{\{i\}}\Big|,
$$
and we say that $F$ is finely-$\mcC^k_b$ if $\big|F\big|_{\mcC^k}$ is finite. 

\ssk

\noindent In example \eqref{EqExample1} above, we have
$$
V'\Big(y_{[0,a]};\dot y_{[0,a]}\Big) = \sum_{i=1}^n \big(1-r_i'(a)\big)\partial_if\big(y_{a-r_1(a)\vee 0},\dots,y_{a-r_n(a)\vee 0}\big)\big(\dot y_{a-r_i(a)}\big),
$$
where $\partial_i$ stands for the differential with respect to the $i^{\textrm{th}}$ coordinate, while in example \eqref{EqExample2} we have
$$
V'\Big(y_{[0,a]};\dot y_{[0,a]}\Big) = \frac{1}{a}\Big(g(y_a)-V\big(y_{[0,a]}\big)\Big),
$$
for $a>0$ and $V'\big(y_{[0,0]};\dot y_{[0,0]}\big) = \big(D_{y_0}g\big)(\dot y_0)$. The path-dependent vector field 
$$
V\Big(y_{[0,a]}\Big) = \int_0^a g(y_r)dr
$$
is also finely-$\mcC^1$. 

\ssk

\noindent  Remark that the restrictions to $\mcC^1\big([0,0],\RR^d\big)\simeq \RR^d\times\RR^d$, of the above two notions of differentiation coincide, as we have 
\begin{equation}
\label{EqV'DV}
V'\big(y;W(y)\big) = (WV)(y)
\end{equation}
for any $y\in\RR^d$ and $W(y)\in\RR^d$.  They do not commute in general: $(DF)'\neq D(F')$. Also, for a map $F$ which is both Fr\'echet differentiable and finely differentiable, with $F'$ Fr\'echet differentiable, we have for any $\mcC^1$ paths $y_\bullet,z_\bullet$ in $\RR^d$ the relation
$$
D_{y_{[0.a]}}F\Big(z_{[0,a]}\Big) = \Big(D_{y_{[0,a]}}\big(F'\big)\Big)\Big(z_{[0,a]} ; \dot y_{[0,a]}\Big) + F'\Big(y_{[0,a]} ; \dot z_{[0,a]}\Big).
$$
Note that if an $\RR^d$-valued path-dependent vector field $V$ over $\RR^d$ is finely-$\mcC^1$ and $f : \RR^d\rightarrow\RR$ is of class $\mcC^2$, then the map $Vf : \mcP \rightarrow \RR$, defined by the formula
$$
(Vf)\Big(y_{[0,a]}\Big) = (D_{y_a}f\big)\Big(V\Big(y_{[0,a]}\Big)\Big)
$$
is finely-$\mcC^1$, and we have
\begin{equation}
\label{EqFormulaVf'}
(Vf)'\Big(y_{[0,a]};\dot y_{[0,a]}\Big) = \Big(D^2_{y_a}f\Big)\Big(\dot y_a,V\Big(y_{[0,a]}\Big)\Big) + (D_{y_a}f)\Big(V'\Big(y_{[0,a]};\dot y_{[0,a]}\Big)\Big).
\end{equation}
for any $y_\bullet\in\mcP^1$. Note also that if $V$ is $\mcC^1$ and $W$ is another finely-$\mcC^1$ path-dependent vector field over $\RR^d$, then
\begin{equation}
\label{EqVWPrime}
(VW)'\Big(y_{[0,a]}\,;\cdot\Big) = (D W)'\Big(y_{[0,a]}\,;\cdot\Big) \left(\overline{V\Big(y_{[0,a]}\Big)}_{[0,a]}\right) + \Big(D_{y_{[0,a]}}V\Big) \left(\overline{V'\Big(y_{[0,a]};\cdot\Big)}_{[0,a]}\right).
\end{equation}

\bigskip

\subsection{Path-dependent rough differential equations}
\label{SubsectionPathDpdtRDEs}

Let $V,V_1,\dots,V_\ell$ be paths dependent vector fields over $\RR^d$, and ${\bfX} = (X,\bbX)$ be a Holder weak geometric $p$-rough path, with $2\leq p<3$, defined on the time interval $[0,T]$. Let $h : [0,T]\rightarrow\RR$ be an $\alpha$-Lipschitz real-valued path, with $\frac{1}{p}+\alpha>1$. We make the following regularity assumptions on the path-dependent vector fields $V,V_i$. 

\medskip

\noindent {\bf (H1) Regularity assumptions}
\begin{itemize}
   \item[-] The drift vector field $V$ is stronly Lipschitz and $\mcC^2_b$.
   \item[-] The $V_i$'s are strongly Lipschitz and have finite $\llparenthesis \cdot\rrparenthesis$-norm, as defined by the formula
$$
\llparenthesis V_i\rrparenthesis := \sum_{k=0}^2 \left\{ \Big|V_i^{(k)}\Big|_{\mcC^{3-k}} + \left\|\Big(V_i^{(k)}\Big)^{\{2-k\}}\right\|_{\mcC^1}\right\}.
$$
\end{itemize}

\ssk

The path-dependent vector fields introduced in examples \eqref{EqExample1} and \eqref{EqExample2} satisfy these assumptions provided $f$ or $g$ are sufficiently regular, with bounded derivatives. 

\bigskip

Given $0\leq s\leq t\leq T$, define $\mu_{ts} : \RR^d\rightarrow \RR^d$, as the well-defined time $1$ map associated with the ordinary path-dependent differential equation 
\begin{equation}
\label{ODE-RDE}
\dot{y}_u = \big(h_t-h_s\big)V\Big(y_{[0,u]}\Big) + X^i_{ts}V_i\Big(y_{[0,u]}\Big) + \frac{1}{2}\XX^{jk}_{ts}[V_j,V_k]\Big(y_{[0,u]}\Big), \quad 0\leq u\leq 1.
\end{equation}

\noindent As noted above, the identification of the restriction to $\mcC\big([0,0],\RR^d\big)$ of $V$ and $V_i$ with vector fields on $\RR^d$ gives their classical meaning to the expressions $V(x), (Vf)(x)$ or $(V_jV_kf)(x)$ in the proposition \ref{PropFundamentalEstimate} below. Set 
$$
\gamma := \min\left\{\frac{3}{p},\alpha+\frac{1}{p}\right\} > 1.
$$

\begin{prop}
\label{PropFundamentalEstimate}
There exists a positive constant $M$ depending polynomially in $\|V\|_{\mcC^2}$ and the $\llparenthesis V_i\rrparenthesis$'s, such that we have 
\begin{equation}
\label{EqFundamentalEstimate}
\Big\|f\circ\mu_{ts} - \Big\{f + \big(h_t-h_s\big)(Vf) + X^i_{ts}\big(V_if\big) + \XX^{jk}_{ts}\big(V_jV_kf\big)\Big\}\Big\|_\infty \leq M\Big(1+\|{\bfX}\|^3\Big)\|f\|_{\mcC^3}|t-s|^\gamma,
\end{equation}
for any real-valued function $f$ of class $\mcC^3$ on $\RR^d$, and any $0\leq s\leq t\leq T$.
\end{prop}

The proof is based on the following elementary inequality
\begin{equation}
\label{EqBasicEstimateYu}
\big|y_u-x\big| \leq m\big(1+\|{\bfX}\|\big)|t-s|^{\frac{1}{p}},
\end{equation}
in which 
$$
m\leq c\left(\|V\|_\infty+\sum_{j=1}^\ell\big(1+\|V_j\|_{\mcC^1}\big)^2\right),
$$ 
for some constant $c$ depending only on $h$. Write $\bbA_{ts}$ for the antisymmetric part of the matrix $\XX_{ts}$. From equation \eqref{ODE-RDE}, we have for any function $f : \RR^d\rightarrow \RR$, of class $\mcC^3$, the identity
\begin{equation}
\label{EqFundamentalIdentity}
\begin{split}
f\big(\mu_{ts}(x)\big) &= f(x) + (t-s) \int_0^1 (Vf)\Big(y_{[0,u]}\Big)du + X^i_{ts}\int_0^1 \Big(V_if\Big)\big(y_{[0,u]}\big)du + \bbA^{jk}_{ts}\int_0^1  \big(V_jV_kf\big)\Big(y_{[0,u]}\Big)du \\
                                  &= f(x) + (t-s)(Vf)(x) + X^i_{ts}\big(V_if\big)(x) + \XX^{jk}_{ts}\big(V_jV_kf\big)(x) + \ep^f_{ts}(x),                                  
\end{split}
\end{equation}
where 
\begin{equation*}
\begin{split}
\ep^f_{ts}(x) &:= \big(h_t-h_s\big)\int_0^1\Big\{(Vf)\Big(y_{[0,u]}\Big)-(Vf)(x)\Big\}du \\
&+ \int_0^1\int_0^u X^i_{ts}X^j_{ts}\Big\{\big(V_if\big)'\Big(y_{[0,r]};V_j(y_\bullet)_{[0,r]}\Big)  - \big(V_if\big)'\Big(y_{[0,0]};V_j\big(y_{[0,0]}\big)\Big)\Big\}\,drdu \\
&+ \int_0^1\int_0^u \Big\{\big(h_t-h_s\big)X^i_{ts}\big(V_if\big)'\Big(y_{[0,r]};V(y_\bullet)_{[0,r]}\Big) + X^i_{ts}\XX^{jk}_{ts}\big(V_if\big)'\Big(y_{[0,r]};\big[V_j,V_k\big](y_\bullet)_{[0,r]}\Big) \Big\}\,drdu
\end{split}
\end{equation*}
We used here the fact that $\bfX$ is weak geometric, so the symmetric part of $\XX_{ts}$ is $\frac{1}{2}X_{ts}\otimes X_{ts}$; formula \eqref{EqFineRegularity} was also applied to $V_if$. 

\medskip

\begin{Dem}
The first term in $\ep^f_{ts}(x)$ is seen to be bounded above by 
$$
c\,\|f\|_{\mcC^1}\,m\|V\|_{\mcC^1}\big(1+\|{\bfX}\|\big)|t-s|^{\alpha+\frac{1}{p}},
$$
using inequality \eqref{EqBasicEstimateYu} and the strong Lipschitz character of $V$. Using \eqref{EqBasicEstimateYu} and the fact that the $V_i'$ are $\mcC^1$, with bounded derivatives, to deal with the second and third terms in $\ep^f_{ts}(x)$, we see that
$$
\Big\|\ep^f_{ts}\Big\|_\infty \leq M\big(1+\|\bfX{}\|^3\big)\|f\|_{\mcC^3}|t-s|^\gamma,
$$
with $M\leq c\sum_{i,j=1}^\ell \Big(1+\big\|V_i'\big\|_{\mcC^1}\Big)\Big(1+\|V_j\|_{\mcC^1}\Big)^3\Big(1+\|V\|_{\mcC^1}\Big)$.
\end{Dem}

\medskip

Note that under the regularity assumptions {\bf (H1)} on the vector fields $V$ and $V_i$'s, the map $\ep^f_{ts}$ is differentiable and we have 
\begin{equation}
\label{EqEstimateDEp}
\Big|D_x\ep^f_{ts}\Big| \leq M'\Big(1+\|{\bfX}\|^4\Big)\|f\|_{\mcC^4}|t-s|^\gamma,
\end{equation}
for all $f\in\mcC^4_b$ and $x\in\RR^d$, where $M'$ depends polynomially on $\|V\|_{\mcC^2}$ and the $\llparenthesis V_i\rrparenthesis$'s.

\begin{thm}
\label{ThmC1ApproximateFlow}
The family of maps $\big(\mu_{ts}\big)_{0\leq s\leq t\leq T}$ is a $\mcC^1$-approximate flow which depends continuously on $\big((s,t);{\bfX},V,V_i\big)$, in the topology of uniform convergence, with respect to the product topology associated with the rough path metric and the metrics $\|\cdot\|_{\mcC^2}$ and $\llparenthesis\cdot\rrparenthesis$.
\end{thm}

\begin{Dem}
The family $\big(\mu_{ts}\big)_{0\leq s\leq t\leq T}$ is made up of functions of class $\mcC^2$, as a direct consequence of classical results on the dependence of solutions to ordinary differential equations \wrt parameters, as recalled above. These results also imply the continuous dependence of $\mu_{ts}$ on $\big((s,t);{\bf X}, V,V_i\big)$. To show that $\big(\mu_{ts}\big)_{0\leq s\leq t\leq T}$ has the {\it $\mcC^1$-approximate flow property} \eqref{EqMuMu}, write, for $0\leq s\leq u\leq t\leq T$,
\begin{equation*}
\mu_{tu}\big(\mu_{us}(x)\big) = \mu_{us}(x) + \big(h_t-h_u\big)V\big(\mu_{us}(x)\big) + X^i_{tu}V_i\big(\mu_{us}(x)\big) + \XX^{jk}_{tu}\big(V_jV_k\big)\big(\mu_{us}(x)\big) + \ep^{\textrm{Id}}_{tu}\big(\mu_{us}(x)\big)
\end{equation*}
and expand each term, so as to get for $\mu_{tu}\big(\mu_{us}(x)\big)$ the expression
\begin{equation*}
\begin{split}
x&+ \big(h_u-h_s\big) V(x) + X^i_{us}V_i(x) + \XX^{jk}_{us}\big(V_jV_k\big)(x) + \ep^\textrm{Id}_{us}(x) \\
  &+ \big(h_t-h_u\big)V(x) + \big(h_t-h_u\big)\Big(V\big(\mu_{us}(x)\big)-V(x)\Big) \\
  &+ X^i_{tu}\Big\{V_i(x) + \big(h_u-h_s\big)\big(VV_i\big)(x) + X^{i'}_{tu}\big(V_{i'}V_i\big)(x) + \XX^{jk}_{tu}\big(V_jV_kV_i\big)(x) + \ep^{V_i}_{us}(x)\Big\} \\
  &+ \XX^{jk}_{tu}\big(V_jV_k\big)(x) + \XX^{jk}_{tu}\Big(\big(V_jV_k\big)\big(\mu_{us}(x)\big) - \XX^{jk}_{tu}\big(V_jV_k\big)(x)\Big).
\end{split}
\end{equation*}
Pairing terms together we see that
$$
\mu_{tu}\big(\mu_{us}(x)\big) = \mu_{ts}(x) + \big\{\cdots\big\},
$$
\noindent where 
\begin{equation*}
\begin{split}
\big\{\cdots\big\} &:= \ep^\textrm{Id}_{ts}(x) + \ep^\textrm{Id}_{us}(x) + \big(h_t-h_u\big)\Big(V\big(\mu_{us}(x)\big)-V(x)\Big) + X^i_{tu}\ep^{V_i}_{us}(x) \\
&\;\;\;+  X^i_{tu}\XX^{jk}_{tu}\big(V_jV_kV_i\big)(x) + \XX^{jk}_{tu}\Big\{\big(V_jV_k\big)\big(\mu_{us}(x)\big) - \XX^{jk}_{tu}\big(V_jV_k\big)(x)\Big\}
\end{split}
\end{equation*}
So we have 
$$
\Big\|\mu_{tu}\big(\mu_{us}(x)\big)-\mu_{ts}(x)\Big\|_\infty \leq M\Big(1+\|{\bfX}\|^3\Big)|t-s|^\gamma
$$
as a consequence of \eqref{EqFundamentalEstimate}, \eqref{EqBasicEstimateYu} and the Lipschitz character of $V$, for some constant $M$ which depends polynomially on $\|V\|_{\mcC^2}$ and the $\llparenthesis V_i\rrparenthesis$'s. To estimate 
$$
D_x\mu_{tu}\circ\mu_{us} - D_x\mu_{ts}
$$
we use \eqref{EqEstimateDEp} for the two terms of the form $\ep_{ba}^\textrm{Id}(x)$, and the regularity assumptions on $V $and the $V_i$'s to get some obvious estimates for
$$
\big(h_t-h_u\big)\Big\{V\big(\mu_{us}(x)\big)-V(x)\Big\}  + X^i_{tu}\XX^{jk}_{tu}\big(V_jV_kV_i\big)(x) + \XX^{jk}_{tu}\Big\{\big(V_jV_k\big)\big(\mu_{us}(x)\big) - \XX^{jk}_{tu}\big(V_jV_k\big)(x)\Big\}.
$$
The terms $X^i_{tu}D_x\ep^{V_i}_{us}$ are dealt with using formula \eqref{EqVWPrime} and the regularity assumptions on the $V_i$'s, which were taylor made to handle that term. This gives an upper bound of the form
$$
\Big|D_x\mu_{tu}\circ\mu_{us} - D_x\mu_{ts}\Big| \leq M\big(1+\|{\bfX}\|^3\big)|t-s|^\gamma,
$$
where the constant $M$ depends polynomially on $\|V\|_{\mcC^2}$ and the $\llparenthesis V_i\rrparenthesis$'s. 

So far, we have only used the fact that the $V_i'$ are $\mcC^1_b$, rather than $\mcC^2_b$, as assumed. This stronger requirement is made to ensure that the decomposition \eqref{EqRegularityBounds} for $D_x\mu_{ts}$ holds, as given by identity \eqref{EqFundamentalIdentity} applied to the identity as function $f$; it follows that the family $\big(\mu_{ts}\big)_{0\leq s\leq t\leq T}$ is indeed a $\mcC^1$-approximate flow on $\RR^d$.
\end{Dem}

\ssk

Write F for the collection $\big(V_1,\dots,V_\ell\big)$.

\begin{defn}
\label{DefnSolRDE}
A \textbf{\emph{flow}} $(\varphi_{ts}\,;\,0\leq s\leq t\leq T)$ on $\RR^d$ is said to \textbf{\emph{solve the path-dependent rough differential equation}}
\begin{equation}
\label{RDE}
d\varphi = Vh(dt) + {\emph F}{\bfX}(dt)
\end{equation}
if there exists a constant $\gamma>1$ independent of $\bfX$ and two possibly $\bfX$-dependent positive constants $\delta$ and $c$ such that
\begin{equation}
\label{EqDefnSolRDE}
\|\varphi_{ts}-\mu_{ts}\|_\infty \leq c\,|t-s|^\gamma
\end{equation}
holds for all $0\leq s\leq t\leq T$ with $t-s\leq\delta$.
\end{defn}

This notion of solution flow follows Davie's definition of a solution to a classical rough differential equation \cite{Davie}, given in terms of its Euler expansion. Proposition \ref{PropFundamentalEstimate} justifies this definition. The following well-posedness result comes as a direct consequence of theorems \ref{ThmConstructingFlows} and \ref{ThmC1ApproximateFlow}.

\begin{thm}[Well-posedness]
\label{ThmMain}
Assume the path-dependent vector fields $V, V_i$ satisfy the regularity assumptions {\bf  (H1)}. Then the path-dependent rough differential equation \eqref{RDE} has a unique solution flow; it depends continuously on $\big((s,t);{\bfX},V,V_i\big)$, in the topology of uniform convergence, with respect to the product topology associated with the rough path metric and the metrics $\|\cdot\|_{\mcC^2}$ and $\llparenthesis\cdot\rrparenthesis$.
\end{thm}

\begin{Dem}
\noindent Note that any solution flow to the rough differential equation \eqref{RDE} depends, by definition, continuously on $(s,t)$ in the topology of uniform convergence since $\mu_{ts}$ does. Use the notations $c_1 \delta$ of section \ref{SectionApproximateFlows} for the constants appearing in the definition of a $\mcC^1$-approximate flow and theorem \ref{ThmConstructingFlows} on their associated flows. It follows from the above estimates and the computations of section 2 in \cite{RPDrivenFlows} giving $\delta$ in terms of $c_1$ and  that we can choose 
$$
c_1=M\Big(1+\|{\bfX}\|^3\Big),\quad \delta = c\,c_1^p,
$$
where the constant $M$ depends polynomially on $\|V\|_{\mcC^2}$ and the $\llparenthesis V_i\rrparenthesis$'s, so we have 
\begin{equation}
\label{EqEstimateVarphiMu}
\|\varphi_{ts}-\mu_{\pi_{ts}}\|_\infty \leq cM^2\Big(1+\|{\bfX}\|^6\Big) T\,\big|\pi_{ts}\big|^\frac{1}{p},
\end{equation}
for any partition $\pi_{ts}$ of $(s,t)\subset [0,T]$ with mesh $\big|\pi_{ts}\big|\leq \delta$, as a consequence of inequality \eqref{EqApproxVarphiMu}. These bounds are uniform in $(s,t)$, for $\bfX$ in a bounded set of the space of weak geometric $p$-rough paths and $V,V_i$'s in $\|\cdot\|_{\mcC^2} / \llparenthesis\cdot\rrparenthesis$-bounded sets. As each $\mu_{\pi_{ts}}$ is a continuous function of $\big((s,t);\bfX,V,V_i\big)$ by theorem \ref{ThmC1ApproximateFlow}, the flow $\varphi$ depends continuously on $\big((s,t);\bfX,V,V_i\big)$.
\end{Dem}

\ssk

\begin{cor}[Vanishing delay]
\label{CorVanishingDelay}
Suppose $V^\lambda$ converges to $V$ in $\mcC^2_b$ and the $V_i^\lambda$'s converge to $V_i$ in $\llparenthesis\cdot\rrparenthesis$-norm, as $\lambda$ tends to $\infty$. If $V\big(y_{[0,u]}\big)$and $V_i\big(y_{[0,u]}\big)$ depend only on $y_u$, for all $0\leq u\leq T$, then the solution flow $\varphi^\lambda$ to the path-dependent rough differential equation 
$$
d\varphi^\lambda = V^\lambda h(dt) + \textrm{\emph{F}}^\lambda{\bfX}(dt)
$$
converges uniformly to the solution flow $\varphi$ of the classical, "Markovian", rough differential equation on flows
$$
d\varphi = Vh(dt) + \textrm{\emph{F}}{\bfX}(dt).
$$
\end{cor}

This corollary applies in particular to path-dependent vector fields of the form \eqref{EqExample1}, with delays $r_1,\dots,r_n$ uniformly bounded above by a constant $o_\lambda(1)$, providing a generalization to our setting of the results of Ferrante and Rovira \cite{FerranteRovira} on stochastic differential equations with delay driven by a fractional Brownian motion.

\begin{rems}
\label{Remarks}
\begin{enumerate}
   \item The setting developped in this work is different from the classical setting of stochastic functional differential equations, as developped for instance in the works of Mohammed \cite{Mohammed}. The main difference is that by taking initial conditions of the form $x_{[0,0]}\in\RR^d$ rather than paths $x_{[-r,0]}\in\mcC\big([0,r],\RR^d\big)$, we can work in $\RR^d$ rather than in the infinite dimensional Banach space $\mcC\big([0,r],\RR^d\big)$. This somehow corresponds to considering here the initial condition as modelling an object without history in the physical system under study. \vspace{0.1cm}  
   
   \item {\bf Stochastic path-dependent differential equations.} We are interested, in this classical setting, in solving the path-dependent stochastic differential equation in $\RR^d$, in Stratonovich form,
   \begin{equation}
   \label{EqStrato}
   dx_t = V\Big(x_{[0,t]}\Big)dt + V_i\Big(x_{[0,t]}\Big){\circ dW^i_t},   
   \end{equation}
   where $W$ is a Brownian motion defined on some probability space $\big(\Omega,\mcF,\big(\mcF_t\big)_{0\leq t\leq T},\PP\big)$; the notation $\big(\mcF_t\big)_{0\leq t\leq T}$ stands here for the completed filtration generated by $W$. This equation has a unique strong solution started from any given point $x_0\in\RR^d$, under sufficient regularity assumptions on $V$ and the $V_i$'s; see for instance \cite{Bichteler}. Given any $\mcC^2_b$ real-valued function $f$ on $\RR^d$, it turns the process
   $$
   M^x_t := f\big(x_t\big)-f\big(x_0\big) - \int_0^t\Big(Vf+\frac{1}{2}V_iV_if\Big)\Big(x_{[0,r]}\Big)dr, \quad\quad   0\leq t\leq T,
   $$
   into a martingale. The latter property characterizes uniquely the solution process $x_\bullet$.
   
   Assume now that $V$ and the $V_i$ satisfy the regularity assumptions {\bf (H1)}, and let ${\bfX} = (W,\bbW)$ denote the Brownian rough path above $W$, and $\varphi$ stand for the solution flow to the path-dependent rough differential equation 
   \begin{equation}
   \label{EqRDEdt}
   d\varphi = Vdt + \textrm{\emph{F}}{\bfX}(dt),
   \end{equation}
   as defined above. Setting $y_t = \varphi_{t0}\big(x_0\big)$, it follows from \eqref{EqFundamentalEstimate} and \eqref{EqDefnSolRDE} that we have
   \begin{equation}
   \label{EqTaylorExp}   
   \big|\Delta^y_{ts}\big| := \Big|f(y_t)-f(y_s) -\Big\{(t-s)(Vf)\Big(y_{[0,s]}\Big) + W^i_{ts}\big(V_if\big)\Big(y_{[0,s]}\Big) + \bbW^{jk}_{ts}\big(V_jV_kf\big)\Big(y_{[0,s]}\Big) \Big\}\Big| \leq c({\bfX})\,|t-s|^\gamma,
   \end{equation}
   for some constant $c(\bfX)$ depending polynomially on the H\"older norm of the rough path $\bfX$; so it defines an integrable random variable. Set $s_i = s+\frac{i}{n}(t-s)$, for $0\leq i\leq n-1$ and $n\geq 2$. As $\sum_{i=0}^{n-1} \Delta^y_{s_{i+1}s_i}$ converges in $\LL^1$ to $M^y_{ts}-\int_s^t\big(V_if\big)\Big(y_{[0,r]}\Big)dr$, as $n$ goes to $\infty$, it follows from the upper bound \eqref{EqTaylorExp} that 
   $$
   M^y_{ts} = \int_s^t\big(V_if\big)\Big(y_{[0,r]}\Big)dW^i_r,
   $$
   is indeed a martingale, so the path $y_\bullet$ constructed from the solution flow to the rough differential equation \eqref{EqRDEdt} coincides with the classical strong solution to the Stratonovich equation \eqref{EqStrato}. 
   
   As illustrated in the recent work \cite{Shevchenko} of Shevchenko, one can solve uniquely some path-dependent It\^o stochastic differential equation of the form 
   $$
   dx_t = a\Big(t,x_{[0,t]}\Big)dt + c\Big(t,x_{[0,t]}\Big)dZ_t + b\Big(t,x_{[0,t]}\Big)dW_t,
   $$
   where $W$ is a Brownian motion and $Z$ an $\alpha$-H\"older process, with $\alpha>\frac{1}{2}$, under weaker assumptions on the path-dependent vector fields than assumptions {\bf (H1)}. Although one can still get some control on the solution process in this setting, one looses the crucial continuity property of the solution map associating $x_\bullet$ to $W_\bullet$. On the other hand, the continuity of the Ito-Lyons map
   $$
   {\bfX}\mapsto \varphi,
   $$
   given by theorem \ref{ThmMain}, provides for free a support theorem for the distribution of the solution $\varphi$ to equation \eqref{RDE}, and a large deviation principle for its "small noise" counterpart
   $$
   \varphi^\epsilon = Vh(dt) + \textrm{\emph{F}}{\bfX}^\epsilon(dt)   
   $$
    where ${\bfX}^\epsilon$ is the rough path associated with $\epsilon W$, as direct consequences of the support theorem for $\bfX$ and the large deviation principle for ${\bfX}^\epsilon$ in rough path topology. See \cite{MohammedZhang} for earlier results in that directions, and the book \cite{FVBook} for an account of these point for classical rough differential equations. Note however that strong regularity assumptions on the vector fields $V,V_i$ are needed to get the continuity of the Ito-Lyons map. 
    
    The pathwise nature of our approach to stochastic path-dependent differential equations, as developped in this section, makes it possible to deal with non-adapted random path-dependent vector fields and initial conditions, which is hard to do in the classical Ito-Skorokhod setting.   \vspace{0.1cm}  
    
   \item {\bf Convergence of Milstein schemes.} Estimate \eqref{EqEstimateVarphiMu} provides a sharp pathwise rate of convergence for the Milstein scheme associated with \eqref{EqStrato}, given here under the form of the maps $\mu_{\pi_{ts}}$. This generalizes to the rough path setting the results of \cite{HuMohammedYan} and \cite{KloedenShardlow} obtained for fractional Brownian motion with ad hoc tools. However, this does not give any insight on the rates of convergence of the associated Euler schemes, as investigated for instance in \cite{GyongySabanis} or \cite{KumarSabanis}. \vspace{0.1cm}
  
   \item Neuenkirch, Nourdin and Tindel obtained in \cite{NeuNourdinTindel} a well-posedness result for a particular case of stochastic delayed differential equation driven by a fractional Brownian motion with Hurst index greater than $\frac{1}{3}$, with no drift and vector fields $V_i$ of the form of example \eqref{EqExample1}. They used for that purpose a variant of Gubinelli's controlled rough paths which requires the introduction of some delayed L\'evy areas. The above results show that the extended notion of rough path introduced in \cite{NeuNourdinTindel} is not needed to get the above general results. Note however that they consider delayed equations with bounded memory and initial condition a path $y_{[-r,0]}\in\mcC\big([-r,0],\RR^d\big)$, as in the classical setting. This is different from what we are doing above.  \vspace{0.1cm}
   
   \item By enlarging the state space to $\RR^d\times\Big(\RR^d\times\big(\RR^d\big)^{\otimes 2}\Big)$, one can deal with path-dependent rough differential equations where the vector fields $V,V_i$ depend not only on $y_{[0,u]}$ but also on ${\bfX}_{[0,u]}$. In so far as $\bfX$ satisfies a linear rough differential equation in $(\RR^d\times\big(\RR^d\big)^{\otimes 2}$ driven by $\bfX$ itself, this requires the use of local $\mcC^1$-approximate flows, as introduced in \cite{RPDrivenFlows}.
   \end{enumerate}
\end{rems}

\section{Path-dependent rough differential equations driven by general H\"older weak geometric $p$-rough paths}
\label{SectionGeneralRoughPaths}

Following \cite{RPDrivenFlows}, we show in this section how the above results of section \ref{SectionPathDependentRDE} on path-dependent rough differential equations driven by a (H\"older weak geometric) $p$-rough path, with $2<p<3$, can easily be adapted to deal with path-dependent rough differential equations driven by any H\"older weak geometric $p$-rough path, for any $2<p<\infty$. 

\medskip

We work in this section with path-dependent vector fields $V,V_1,\dots,V_\ell$ over $\RR^d$, satisfying the following regularity assumptions.

\medskip

\noindent {\bf (H2) Regularity assumptions.}
\begin{itemize}
   \item[-] The path-dependent drift vector field $V$ is strongly Lipschitz and $\mcC^2_b$.
   \item[-] The $V_i$'s are strongly Lipschitz and have finite $\llparenthesis \cdot\rrparenthesis$-norm, as defined by the formula
$$
\llparenthesis V_i\rrparenthesis := \sum_{k=0}^{[p]} \left\{ \Big|V_i^{(k)}\Big|_{\mcC^{[p]+1-k}} + \left\|\Big(V_i^{(k)}\Big)^{\{[p]-k\}}\right\|_{\mcC^1}\right\}.
$$
\end{itemize}

\subsection{Constructing a $\mcC^1$-approximate flow.} 
\label{SubsectionGeneralC1ApproximateFlow}

Let $h : [0,T]\rightarrow\RR$ be an $\alpha$-Lipschitz real-valued path, with $\frac{1}{p}+\alpha>1$. Under the above regularity assumptions, and given a tuple $I=(i_1,\dots,i_r)$, with $2\leq n\leq [p]$, we define a vector field $V_{[I]}$ and a path-dependent differential operator $V_I $ on $\RR^d$ setting
$$
V_{[I]} := \Big[V_{i_1},\big[V_{i_2},\dots[V_{i_{r-1}},V_{i_r}]\big]\dots\Big], \quad \textrm{and }\quad  V_I = V_{i_1}V_{i_2}\cdots V_{i_r}.
$$

\ssk

 Given $0\leq s\leq t\leq T$, let $\mu_{ts}$ be the well-defined time $1$ map associated with the path-dependent ordinary differential equation 
\begin{equation}
\label{ODE-RDEGeneral}
\dot y_u = \big(h_t-h_s\big)V\Big(y_{[0,u]}\Big) + \sum_{|I|\leq [p]} \Lambda_{ts}^I V_{[I]}\Big(y_{[0,u]}\Big), \quad\quad    0\leq u\leq 1.
\end{equation}
Note that this equation reduces to equation \eqref{ODE-RDE} of section \ref{SectionPathDependentRDE}, when $2<p<3$. Set $\gamma=\min\left\{\frac{[p]+1}{p},\frac{1}{p}+\alpha\right\}>1$.
 
\begin{prop}
There exists a constant $M$ depending polynomially on $\|V\|_{\mcC^2}$ and the $\llparenthesis V_i\rrparenthesis$'s, such that the inequality
\begin{equation}
\label{EqFundEstimateGeneral}
\left\|f\circ\mu_{ts} - \left\{f + (t-s)Vf + \sum_{r=1}^{[p]} \sum_{I\in\llbracket 1,\ell\rrbracket^r} X^I_{ts}V_If\right\}\right\|_\infty \leq M\big(1+\|\bfX\|^\gamma\big)\,\|f\|_{\mcC^{[p]+1}}\,|t-s|^\gamma,
\end{equation}
 holds for any real-valued function $f\in\mcC^{[p]+1}$, and any $0\leq s\leq t\leq T$. 
\end{prop}

The proof of this proposition and the next one are based on the elementary identity \eqref{EqExactFormulaMuTs} below, obtained by applying repeatedly the identity
$$
f(y_r) = f(x) + \big(h_t-h_s\big)\int_0^1(Vf)\Big(y_{[0,u]}\Big)\,du + \sum_{1\leq r\leq [p]}\sum_{I\in\llbracket 1,\ell\rrbracket^r}\Lambda^I_{ts}\int_0^1\big(V_{[I]}f\big)\Big(y_{[0,u]}\Big)\,du,
$$
and by separating the terms according to their size in $|t-s|$. Given some tuples $I_1,\dots,I_k$, with $2\leq k\leq [p]$, and a real-valued function $f$ of class $\mcC^{[p]+1}$, define inductively 
\begin{equation*}
\begin{split}
\Big(V_{[I_{k+1}]}\ast V_{[I_k]}\ast\,\cdots\,\ast V_{[I_1]}f\Big)\big(z_{[0,u]}\big) &:= V_{[I_{k+1}]}\ast\Big(V_{[I_k]}\ast\,\cdots\,\ast V_{[I_1]}f\Big)\big(z_{[0,u]}\big) \\
&= \Big(V_{[I_k]}\ast\,\cdots\,\ast V_{[I_1]}f\Big)'\Big(z_{[0,u]},V_{[I_{k+1}]}(z_\bullet)_{[0,u]}\Big),
\end{split}
\end{equation*}
for any $\mcC^1$ path $z_{[0,u]}$, starting with formula \eqref{EqFormulaVf'}. Since the two differentiation operations $D$ and $'$ on path-dependent vector fields coincide on $\mcC^1\big([0,0],\RR^d\big)\simeq\RR^d\times\RR^d$, as said above in \eqref{EqV'DV}, we have
$$
\Big(V_{[I_k]}\ast\,\cdots\,\ast V_{[I_1]}f\Big)\big(z_{[0,0]}\big) = \big(V_{[I_k]}\cdots V_{[I_1]}f\big)\big(z_0\big),
$$
where the $V_{[I]}$'s are identified in the right hand side with vector fields on $\RR^d$, and $z_{[0,0]}$ is identified with $z_0$. Set $\Delta_n := \big\{(s_1,\dots,s_n)\in [0,T]^n\,;\,s_1\leq \cdots\leq s_n\big\}$, for $2\leq n\leq [p]$. For a function $f$ of class $\mcC^{[p]+1}$,  we have
\begin{equation}
\label{EqExactFormulaMuTs}
\begin{split}
f&\big(\mu_{ts}(x)\big) = f(x) + \big(h_t-h_s\big)\big(Vf\big)(x) + \sum_{k=1}^n \frac{1}{k!}\sum_{|I_1|+\dots+|I_k|\leq [p]} \left(\prod_{m=1}^k \Lambda^{I_m}_{ts}\right) \big(V_{[I_k]}\cdots V_{[I_1]}f\big)(x) \\
&+\sum_{|I_1|+\dots+|I_n|\leq [p]} \left(\prod_{m=1}^n \Lambda^{I_m}_{ts}\right) \int_0^1 \Big\{\big(V_{[I_n]}\ast\cdots \ast V_{[I_1]}f\big)\Big(y_{[0,s_n]}\Big)-\big(V_{[I_n]}\ast\cdots \ast V_{[I_1]}f\big)(x)\Big\}{\bf 1}_{\Delta_n}\,ds_n\dots ds_1 \\
&+ \big(h_t-h_s\big)\int_0^1 \big\{\big(Vf\big)\Big(y_{[0,r]}\Big)-\big(Vf\big)(x)\big\}dr \\
&+ \sum_{k=1}^{n-1}\sum_{|I_1|+\dots+|I_k|\geq [p]+1} \left(\prod_{m=1}^k \Lambda^{I_m}_{ts}\right) \big(V_{[I_k]}\cdots V_{[I_1]}f\big)(x) \\
&+ \big(h_t-h_s\big)\sum_{k=1}^{n-1}\sum_{I_1,\dots,I_k} \left(\prod_{m=1}^k \Lambda^{I_m}_{ts}\right) \int_0^1 \big(V\ast V_{[I_k]}\ast\cdots \ast V_{[I_1]}f\big)\Big(y_{[0,s_k]}\Big){\bf 1}_{\Delta_k}\,ds_k\dots ds_1 \\
&+ \sum_{|I_1|+\dots+|I_n|\geq [p]+1} \left(\prod_{m=1}^k \Lambda^{I_m}_{ts}\right) \int_0^1 \Big\{\big(V_{[I_n]}\ast\cdots \ast V_{[I_1]}f\big)\Big(y_{[0,s_n]}\Big)-\big(V_{[I_n]}\ast\cdots \ast V_{[I_1]}f\big)(x)\Big\}{\bf 1}_{\Delta_n}\,ds_n\dots ds_1
\end{split}
\end{equation}
We denote by $\ep^{f\,;\,n}_{ts}(x)$ the sum of the last four lines, made up of terms of size at least $|t-s|^\gamma$. In the case where $n=[p]$, the terms in the second line involve only some tuples $I_j$ with $|I_j|=1$, so the elementary estimate 
\begin{equation}
\label{EqElementaryEstimateYu}
\big|y_r-x\big| \leq m\big(1+\|{\bfX}\|^{[p]+1}\big)|t-s|^{1/p}, \quad 0\leq r\leq 1,
\end{equation}
can be used to control the increment in the integral, showing that this second line is of order $|t-s|^\gamma$; we include it in the remainder $\ep^{f\,;\,[p]}_{ts}(x)$. As in the proof of proposition \ref{PropFundamentalEstimate}, it is elementary to see that the constant $m$ depends polynomially on $\|V\|_{\mcC^2}$ and the $\llparenthesis V_i\rrparenthesis$'s.

\bigskip

\begin{Dem}
Applying the above formula for $n=[p]$, together with the fact that $\exp({\bf \Lambda}) = \bfX$, we get the identity
\begin{equation*}
f\big(\mu_{ts}(x)\big) = f(x) + \big(h_t-h_s\big)\big(Vf\big)(x) + \sum_{I} X^I_{ts}\big(V_If\big)(x) + \ep^{f\,;\,[p]}_{ts}(x).
\end{equation*}
It is clear on the formula for $\ep^{f\,;\,[p]}_{ts}(x)$ that its absolute value is bounded above by a constant multiple of $\big(1+\|{\bfX}\|^\gamma\big)|t-s|^\frac{\gamma}{p}$, for a constant depending polynomially on $\|V\|_{\mcC^2}$ and $\llparenthesis V_i\rrparenthesis$, and $f$ as in \eqref{EqFundEstimateGeneral}. 
\end{Dem}

\bigskip

Note that under the regularity assumptions {\bf (H2)} on the vector fields $V$ and $V_i$'s, the map $\ep^{f;[p]}_{ts}$ is differentiable and we have 
\begin{equation}
\label{EqEstimateDEp}
\Big|D_x\ep^f_{ts}\Big| \leq M'\Big(1+\|{\bfX}\|^{[p]+1}\Big)\|f\|_{\mcC^{[p]+1}} |t-s|^\gamma,
\end{equation}
for all $f\in\mcC^4_b$ and $x\in\RR^d$, where $M'$ depends polynomially on $\|V\|_{\mcC^2}$ and the $\llparenthesis V_i\rrparenthesis$'s. This is the key remark for proving the next proposition, whose proof is identical to the proof of proposition $11$ in \cite{RPDrivenFlows}. We re-write it here for the reader's convenience.

\begin{thm}
\label{ThmSummaryPropertiesGeneralCase}
The family of maps $\big(\mu_{ts}\big)_{0\leq s\leq t\leq T}$ is a $\mcC^1$-approximate flow which depends continuously on $\big((s,t);{\bfX},V,V_i\big)$, in the topology of uniform convergence, with respect to the product topology associated with the rough path metric and the metrics $\|\cdot\|_{\mcC^2}$ and $\llparenthesis \cdot\rrparenthesis$.
\end{thm}

\begin{Dem}
As in the proof of theorem \ref{ThmMain}, it is elementary to see that $\mu$ depends continuously on $\big((s,t);{\bfX},V,V_i\big)$. To see that it is a $\mcC^1$-approximate flow, we first use formula \eqref{EqExactFormulaMuTs} to write
\begin{equation}
\mu_{tu}\big(\mu_{us}(x)\big) = \mu_{us}(x) + \big(h_t-h_u\big)V\big(\mu_{us}(x)\big) + \sum_I X^I_{tu} V_I\big(\mu_{us}(x)\big) +\ep^{\textrm{Id}\,;\,[p]}_{tu}\big(\mu_{us}(x)\big).
\end{equation}
We deal with the term $\big(h_t-h_u\big)V\big(\mu_{us}(x)\big)$ using \eqref{EqEstimateDEp} and the Lipschitz character of $V$. The remainder $\ep^{\textrm{Id}\,;\,[p]}_{tu}\big(\mu_{us}(x)\big)$ has, under the regularity assumptions {\bf (H2)} a $\mcC^1$-norm bounded above by $M\big(1+\|{\bfX}\|^\gamma\big)^2|t-u|^\gamma$, for some constant $M$ depending polynomially on  $\|V\|_{\mcC^2}$ and the $\llparenthesis V_i\rrparenthesis$'s.

\ssk

To deal with the terms $X^I_{tu} V_I\big(\mu_{us}(x)\big)$, with $|I|=k$, we use formula \eqref{EqExactFormulaMuTs} with $n=[p]-k$, to develop $V_I\big(\mu_{us}(x)\big)$. We have
\begin{equation*}
\begin{split}
&V_I\big(\mu_{us}(x)\big) = V_I(x) + \big(h_u-h_s\big)\big(VV_I\big)(x) + \sum_{j=1}^{[p]-k} \frac{1}{j!}\sum_{|I_1|+\dots+|I_j|\leq [p]} \left(\prod_{m=1}^j \Lambda^{I_m}_{us}\right) \big(V_{[I_j]}\cdots V_{[I_1]}V_I\big)(x)  + \ep^{V_I\,;\,{p}-k}_{us}(x) \\
&+\sum_{|I_1|+\dots+|I_{[p]-k}|\leq [p]} \left(\prod_{m=1}^{[p]-k} \Lambda^{I_m}_{us}\right) \int_0^1 \Big\{\Big(V_{[I_{[p]-k}]}\cdots V_{[I_1]}V_I\Big)\big(y_{s_{[p]-k}}\big)-\Big(V_{[I_{[p]-k}]}\cdots V_{[I_1]}V_I\Big)(x)\Big\}{\bf 1}_{\Delta_{[p]-k}}\,ds_{[p]-k}\dots ds_1.
\end{split}
\end{equation*}
Using \eqref{EqElementaryEstimateYu}, we see that $X^I_{tu}$ times the second line above has a $\mcC^1$-norm bounded above by $M\big(1+\|{\bfX}\|^{[p]+1}\big)^2|t-s|^\gamma$. Writing 
$$
\mu_{us}(x) = x + \big(h_u-h_s\big)V(x) + \sum_{I} X^I_{us}V_I(x) + \ep^{\textrm{Id}\,;\,[p]}_{us}(x),
$$
it is then straightforward to use the identities $\exp {\bf \Lambda}_{us} =\bfX_{us}$ and ${\bfX}_{ts} = {\bfX}_{us}{\bfX}_{tu}$, to see that
$$
\mu_{tu}\big(\mu_{us}(x)\big) = \mu_{ts}(x) + \ep_{ts}(x),
$$
with a remainder $\ep_{ts}$ with a $\mcC^1$-norm bounded above by $M\big(1+\|{\bfX}\|^\gamma\big)^2|t-s|^\gamma$.

We see that $\mu_{ts}$ is a $\mcC^2_b$ perturbation of the identity using \eqref{EqFundamentalIdentity} applied to the identity as function $f$.
\end{Dem}

\subsection{ Well-posedness result for path-dependent rough differential equations.} 
\label{SubsectionGeneralWellPosedness}

Defining the solution flow to the path-dependent rough differential equation 
$$
d\varphi = Vh(dt) + \textrm{F}{\bfX}(dt)
$$
as in definition \ref{DefnSolRDE}, theorems \ref{ThmConstructingFlows} and \ref{ThmSummaryPropertiesGeneralCase} imply, as in section \ref{SectionPathDependentRDE}, the following well-posedness result, which contains theorem \ref{ThmMain} as a special case. Its proof is  identical to the proof of theorem \ref{ThmMain} given in section \ref{SectionPathDependentRDE}.

\begin{thm}[Well-posedness]
\label{ThmGeneralMain}
Let $\bfX$ be a H\"older weak geometric $p$-rough path, with $2<p<\infty$, and $V, V_i$ be path-dependent vector fields satisfying the regularity assumptions {\bf  (H2)}. The path-dependent rough differential equation 
$$
d\varphi = Vh(dt) + \textrm{\emph{F}}{\bfX}(dt)
$$
has a unique solution flow; it depends continuously on $\big((s,t);{\bfX},V,V_i\big)$, in the topology of uniform convergence, with respect to the product topology associated with the rough path metric and the metrics $\|\cdot\|_{\mcC^2}$ and $\llparenthesis\cdot\rrparenthesis$.
\end{thm}

Corollary \ref{CorVanishingDelay} on path-dependent rough differential equations with vanishing delay and remark \ref{Remarks} (2) on the speed of convergence of the associated Milstein schemes apply here as well.

\begin{rem}
Tindel and Torreccila have pushed in \cite{TindelTorrecilla} their machinery of extended controlled paths, with delayed areas and related multiple integrals, to deal with delayed rough differential driven by fractional Brownian motion with Hurst index greater than $\frac{1}{4}$. The results of that section show that the classical setting of rough paths is sufficient to get more general results, when used with the machinery of $\mcC^1$-approximate flows. \vspace{0.1cm}
\end{rem}


\bigskip
\bigskip

\begin{thebibliography}{10}

\bibitem{Mohammed}
S.A.E. Mohammed.
\newblock Stochastic functional differential equations.
\newblock {\em Research Notes in Math.}, 99, 1984.

\bibitem{Mohammed98}
S.A.E. Mohammed.
\newblock Stochastic differential systems with memory: theory, examples and applications.
\newblock {\em Stochastic Analysis and Related Topics VI}, (Decreusefond, Gjerde, Oksendal and Ust\"unel, eds), 1--77, Birkhauser, Boston, 1998.

\bibitem{ItoNisio}
Ito, K. and Nisio, M.
\newblock On stationary solutions of a stochastic differential equation.
\newblock {\em J. Math. Kyoto Univ.}, 4:1--75, 1964.

\bibitem{Kushner}
Ito, K. and Nisio, M.
\newblock On the stability of processes defined by stochastic differential-difference equations.
\newblock {\em J. Diff. Eq.}, 4:424--443, 1968.

\bibitem{Lyons97}
Lyons, T.
\newblock Differential equations driven by rough signals.
\newblock {\em Rev. Mat. Iberoamericana}, 14 (2):215--310, 1998.

\bibitem{RPDrivenFlows}
I. Bailleul.
\newblock Flows driven by rough paths.
\newblock {\em Submitted}, 2012.

\bibitem{Gubinelli}
Gubinelli, M.
\newblock Controlling rough paths.
\newblock {\em J. Funct. Anal.}, 216(1):86--140, 2004.

\bibitem{FdlP1}
Feyel, D. and de La Pradelle, A.
\newblock Curvilinear integrals along enriched paths.
\newblock {\em Electron. J. Probab.}, 11:860--892, 2006.

\bibitem{LyonsStFlour}
Lyons, T.J. and Caruana, M. and L{\'e}vy, Th.
\newblock Differential equations driven by rough paths.
\newblock {\em Lecture Notes in Mathematics}, 1908, Springer 2007.

\bibitem{FdlP2}
Feyel, D. and de La Pradelle, A. and Mokobodzki, G.
\newblock A non-commutative sewing lemma.
\newblock {\em Electron. Commun. Probab.}, 13:24--34, 2008.

\bibitem{Davie}
Davie, A. M.
\newblock Differential equations driven by rough paths: an approach via discrete approximation.
\newblock {\em Appl. Math. Res. Express. AMRX}, (2), 2007.

\bibitem{FerranteRovira}
M. Ferrante and C. Rovira.
\newblock Convergence of delay differential equations driven by fractional Brownian motion.
\newblock {\em J. Evol. Equ.}, 10: 761--783, 2010.

\bibitem{Bichteler}
Bichteler, K.
\newblock Stochastic integration with jump.
\newblock {\em CUP, Cambridge Studies in Advanced Mathematics}, 120, 2010.

\bibitem{Shevchenko}
G. Shevchenko.
\newblock Mixed stochastic delay equations.
\newblock {\em Preprint}, arXiv:1306.0590 , 2013.

\bibitem{MohammedZhang}
Mohammed, S.E.A and Zhang, T.
\newblock Large deviations for stochastic systems with memory.
\newblock {\em Discrete Contin. Dyn. Syst.} Ser. B6 (4): 881--893, 2006.

\bibitem{FVBook}
Friz, P. and Victoir, N.
\newblock Multidimensional stochastic processes as rough paths.
\newblock {\em CUP, Encyclopedia of Mathematics and its Applications}, 89, 2010.

\bibitem{HuMohammedYan}
Y. Hu and S.A.E. Mohammed and Y. Yan.
\newblock Discrete-time approximations of stochastic delay equations: The Milstein scheme.
\newblock {\em Stoch. An. Prob.}, 32: 265--314, 2004.

\bibitem{KloedenShardlow}
P.E. Kloeden and T. Shardlow.
\newblock The Milstein scheme for stochastic delay differential equations without using anticipative calculus.
\newblock {\em Stoch. Anal. Appl.}, 30(2) : 181--202, 2012.

\bibitem{GyongySabanis}
I. Gy\"ongy and S. Sabanis.
\newblock Strong convergence of Euler approximations of stochastic differential equations with delay under local Lipschitz conditions.
\newblock {\em Arxiv.}, 1212.3567, 2013.

\bibitem{KumarSabanis}
C. Kumar and S. Sabanis.
\newblock A note on Euler approximation for stochastic differential equations with delay.
\newblock {\em Arxiv.}, 1303.0017, 2013.

\bibitem{NeuNourdinTindel}
A. Neuenkirch and I. Nourdin and S. Tindel.
\newblock Delay equations driven by rough paths.
\newblock {\em Elec. Journal Prob.}, 13(67) : 2031--2068, 2008.

\bibitem{TindelTorrecilla}
S. Tindel and Torrecilla.
\newblock Some differential systems driven by fractional Brownian motion with Hurst parameter greater than 1/4.
\newblock {\em To appear in Progress in Probability}, (volume dedicated to A.S. Ustunel), 2009.


\end{thebibliography}
\end{document}